\documentclass[a4paper]{amsart}
\usepackage{epsfig} 
\usepackage{graphicx}
\usepackage{multirow}
\usepackage{amsmath,amssymb,amsfonts}
\usepackage{amsthm}
\usepackage{mathrsfs}
\usepackage[title]{appendix}
\usepackage{xcolor}
\usepackage{textcomp}
\usepackage{manyfoot}
\usepackage{booktabs}
\usepackage{algorithm}
\usepackage{algorithmicx}
\usepackage{algpseudocode}
\usepackage{listings}
\usepackage{hyperref}
\usepackage{setspace}
\usepackage{lineno}
\usepackage[left=.9in, right=.9in, top=.9in, bottom=.9in]{geometry}
\usepackage{subcaption}
\usepackage{float}

\usepackage{enumitem}
\theoremstyle{thmstyleone}
\newtheorem{thm}{Theorem}
\newtheorem{lem}{Lemma}[section]
\newtheorem{prop}[lem]{Proposition} 
\newtheorem{cor}[lem]{Corollary}
\usepackage[numbers
]{natbib}

\theoremstyle{thmstyletwo}

\newtheorem{remark}[lem]{Remark}

\theoremstyle{thmstylethree}
\newtheorem{defn}{Definition}

\numberwithin{equation}{section}

\newcommand{\real}{\mathbb{R}}
\newcommand{\complex}{\mathbb{C}}

\newcommand{\integer}{\mathbb{Z}}
\newcommand{\ddd}{\mathbb{D}}

\newcommand{\prtblen}{perturbation length}
\newcommand{\defeq}{:=}

\DeclareMathOperator{\diag}{diag}

\DeclareMathOperator{\oJ}{\mathcal{J}_L}
\DeclareMathOperator{\defJz}{{\mathcal{J}}_{\hspace{-1pt}0}}
\DeclareMathOperator{\rank}{rank}
\DeclareMathOperator{\spec}{spec}
\DeclareMathOperator{\solve}{solve}
\DeclareMathOperator{\nuL}{L}
\DeclareMathOperator{\nuLg}{\ell}

\begin{document}

\title{Inverse scattering for the one-dimensional discrete Schr\"odinger operator: a computational perspective}
\author{{Peter C.} {Gibson}}
\address{{Department of Mathematics \& Statistics}, {York University}, {{Toronto}, {M3J 1P3}, {ON}, {Canada}}}
\email{pcgibson@yorku.ca}

\author{Basant Lal Sharma}
\address{{Department of Mechanical Engineering}, {Indian Institute of Technology Kanpur}, {{Kanpur}, {208016}, {UP}, {India}}}
\email[Corresponding author]{bls@iitk.ac.in}

\maketitle

\begin{abstract}
{
We revisit inverse scattering for the one-dimensional discrete Schr\"odinger operator with compactly supported perturbation. 
We consider a computational perspective, adapted to the case of a digitally sampled reflection coefficient, to determine the perturbation in the form of a Jacobi matrix, including
a computation of the length of the smallest interval supporting the perturbation.
A new representation of scattering data in terms of orthogonal polynomials underpins the computational results. 
A key technical role is played by polynomials on the circle constructed from pairs of successive orthogonal polynomials on the line. 
The inverse algorithm is demonstrated on numerical examples.
}
\end{abstract}



\section{Introduction\label{sec-introduction}}

In the present article, we consider the one-dimensional discrete Schr\"odinger equation
\begin{equation}\label{discrete}
-\Delta u+vu=\sigma^2u,
\end{equation}
for $u(\cdot,\sigma)\in\ell_\infty(\integer)$,  where $\sigma$ is the spectral parameter, $\Delta$ denotes the one-dimensional discrete Laplacian
\begin{equation}\label{laplacian}
\Delta u(n,\cdot)=u(n-1,\cdot)-2u(n,\cdot)+u(n+1,\cdot)\qquad(n\in\integer),
\end{equation}
$v:\integer\rightarrow\real$ is a compactly supported potential and independent of the spectral parameter $\sigma$. 
Without loss of generality, we suppose that the integer lattice is labelled such that the term $vu$ vanishes outside the set $\{1,\ldots, {\nuL}\}$
for some integer ${\nuL}\geq 1$. 

The operator $-\Delta +v$ in \eqref{discrete} is a 
special case of Jacobi operators
\begin{equation}\label{J}
\oJ u(n,\sigma)=-b_{n-1}u(n-1,\sigma)+a_nu(n,\sigma)-b_nu(n+1,\sigma)\qquad(n\in\integer), 
\end{equation}
for $u(\cdot,\sigma)\in\ell_\infty(\integer)$, $b_n>0$, $a_n\in\real$ for all $n\in\integer$,
with the restriction that
\begin{equation}\label{conditions}
b_n=1\mbox{ unless }1\leq n\leq {\nuL}-1,\quad\mbox{ and }\quad a_n=2 \mbox{ unless }1\leq n\leq {\nuL},
\end{equation}
for some integer ${\nuL}\geq 1$.
We shall focus on the perturbed Schr\"odinger equation
\begin{equation}\label{discrete-2}
{\oJ}u=\sigma^2u.
\end{equation}
Henceforth $\oJ$, with structure as above, will be referred to as a perturbed Schr\"odinger operator. 
It is clear that ${\oJ}$ is completely determined by the symmetric, tri-diagonal $\nuL\times\nuL$ matrix 
\begin{equation}\label{Jmat}
\begin{pmatrix}a_1&-b_1&0&\cdots&0\\
-b_1&a_2&-b_2&\ddots&\vdots\\
0&-b_2&\ddots&\ddots&0\\
\vdots&\ddots&\ddots&a_{\mathrm{L}-1}&-b_{\mathrm{L}-1}\\
0&\cdots&0&-b_{\mathrm{L}-1}&a_{\mathrm{L}}
\end{pmatrix}.
\end{equation}
The least ${\nuL}$ such that (\ref{conditions}) holds will be referred to as the {\em {\prtblen}} for ${\oJ}$.  
Evidently, 
(\ref{discrete-2})
reduces to the one-dimensional discrete Schr\"odinger 
equation 
(\ref{discrete}) 
with compactly supported potential
when $b_n\equiv 1$ and $\nuL\ge1$. 
Note that $\oJ$ is a compactly supported perturbation of \begin{equation}\label{J0}
{\defJz}\defeq-\Delta.
\end{equation}
If $\oJ=\defJz$ in (\ref{discrete-2}), the ansatz $u(n,\sigma)=e^{ikn}$ leads directly to the \emph{dispersion relation}
\begin{equation}\label{dispersion}
2-\sigma^2=e^{ik}+e^{-ik},
\end{equation}
as does the ansatz $u(n,\sigma)=e^{-ikn}$.    
For $(\sigma,k)\in [-2,2]\times[-\pi,\pi]$, the dispersion relation \eqref{dispersion} determines a pair of curves corresponding to the respective constraints $\sigma k\geq 0$ and $\sigma k\leq 0$.  The convention $\sigma k\geq 0$ establishes a bijective correspondence 
$[-2,2]\ni \sigma\leftrightarrow k\in[-\pi,\pi],$
whereby we may regard $k$ and $\sigma$ as functions of one another. In particular, $k$ has the explicit form
\begin{equation}\label{k-of-sigma}
k=k(\sigma)\defeq2\arcsin(\sigma/2)\qquad(-2\leq \sigma\leq 2). 
\end{equation}In the present article, 
\begin{equation}\label{notnk}
\text{whenever $k$ appears it represents the expression $k(\sigma)$.}
\end{equation}
Any solution to (\ref{discrete-2}) of the form 
\begin{equation}\label{right}
u(n,\sigma)\sim\beta(k)e^{+ikn}\quad\mbox{ as }\quad n\rightarrow\infty,
\end{equation}
i.e., which eventually takes the form of a purely right travelling wave, is said to be outgoing to the right, or to have no incoming radiation from the right.  Similarly, a solution of the form 
\begin{equation}\label{left}
u(n,\sigma)\sim \alpha(k)e^{-ikn}\text{ as }n\rightarrow-\infty,
\end{equation}
is said to be outgoing to the left. 

The scattering matrix for (\ref{discrete-2}) is defined in terms of special solutions that are outgoing to the right or to the left, as follows---precisely analogous to the scattering matrix for the continuous one-dimensional Schr\"odinger operator \cite{DeTr:1979}. 
For each $0\neq\sigma\in[-2,2]$, equation (\ref{discrete-2}) has uniquely determined solutions of the form
\begin{align}
f_1(n,\sigma)&\sim\left\{
\begin{array}{cc}
e^{ikn}&\mbox{ as }n\rightarrow+\infty\\
\dfrac{1}{T_2(k)}e^{ikn}+\dfrac{R_2(k)}{T_2(k)}e^{-ikn}&\mbox{ as }n\rightarrow-\infty
\end{array}\right.\label{f1}\\
f_2(n,\sigma)&\sim\left\{
\begin{array}{cc}
\dfrac{R_1(k)}{T_1(k)}e^{ikn}+\dfrac{1}{T_1(k)}e^{-ikn}&\mbox{ as }n\rightarrow+\infty\\
e^{-ikn}&\mbox{ as }n\rightarrow-\infty
\end{array}\right..
\label{f2}
\end{align}
In the present context the asymptotic condition $n\rightarrow\pm\infty$ in the above definitions evidently simplify to $n>L$ or $n<1$ as appropriate. 
The scattering (or $S$) matrix for $\oJ$ is defined to be the $2\times2$ $k$-dependent matrix
\begin{equation}\label{S}
S\defeq\begin{pmatrix}T_1&R_2\\ R_1&T_2\end{pmatrix},
\end{equation}
where $R_2$ is called the left reflection coefficient for ${\oJ}$, and $R_1$ the right reflection coefficient; $T_1$ and $T_2$ turn out to coincide and are called the transmission coefficient.  
\begin{remark}
    $T_2f_1$, for example, is a solution to (\ref{discrete-2}) resulting from a right traveling wave $e^{ikn}$ coming from $-\infty$ toward the interval $\{1,\ldots,{\nuL}\}$ supporting the perturbation of ${\oJ}$ from ${\defJz}$.  The incoming wave partially reflects off the perturbation, resulting in the leftward travelling wave $R_2e^{-ikn}$.  The incoming wave also partially transmits through the perturbation, emerging as the wave $T_2e^{ikn}$ which travels rightward to $+\infty$.  A similar interpretation holds for $T_1f_2$. 
\end{remark}

Importantly,
$f_1$ and $f_2$ attain their asymptotic structure outside the support of the perturbation.  
Indeed, (\ref{J}) and (\ref{conditions}) imply 
\begin{equation}\label{asymptotic}
\begin{split}
f_1(n,\sigma)&=\frac{1}{T_2(k)}e^{ikn}+\frac{R_2(k)}{T_2(k)}e^{-ikn}\quad\mbox{ for every } n\leq 1, \mbox{ and }\\
f_1(n,\sigma)&=e^{ikn}\quad\mbox{ for every } n\geq {\nuL}. 
\end{split}
\end{equation}

We now state precisely the two novel computational problems addressed in the present article. 
\subsection*{Problem statement}
Given data $R_2(k_j)$, sampled on a set of evenly spaced sample points $-\pi<k_j<\pi$, the aim is to: 
\begin{align}
    &\label{prob1}\text{ Compute the {\prtblen} $\nuL$ for $\oJ$.}\tag{Problem 1}\\
    &\label{prob2} \text{ Reconstruct $\oJ$ by a finite procedure.\tag{Problem 2}    }
\end{align}

There is an extremely rich volume of theoretical results for the inverse scattering problems and their generalizations for equation \eqref{J}, with hypothesis more general than those associated with (\ref{discrete-2}); see, for example, 
\cite{Guseinov1976a,
Guseinov1976b,
Guseinov1977,
Atkinson1964,
Hochstadt1974,
Yurko1995,
Yurko1996,
Yurko2026,
guseinov1995inverse,
Papanicolaou2011}.
Existing theory is not readily applicable to the above computational problems however, which require novel methods. While scattering data is known to determine the operator, direct computation of the perturbation length from scattering data (without first determining the operator) is a completely new problem.
In this article, we aim to derive effective methods to compute the perturbation length ${\nuL}$ and the operator $\oJ$ (that is, effectively the matrix $J_{\mathrm{L}}$ \eqref{Jmat}), given a sampled version of the scattering matrix for (\ref{discrete-2}). 
Determination of ${\nuL}$ from scattering data is specific to the context in which $\oJ$ is a compactly supported perturbation of $\defJz$.
The 
inverse scattering problem can be posed as the question of determination of the operator $\oJ$, given the scattering matrix $S$. As is well known, either one of the reflection coefficients $R_1, R_2$, determines $S$. 
In this article, we work with the left reflection coefficient $R_2$ as given data, corresponding to an unknown perturbed Schr\"odinger operator $\oJ$ satisfying (\ref{J}) and (\ref{conditions}).
From a purely mathematical point of view, scattering for (\ref{discrete-2}) is intriguing as a hybrid setting, somewhere between the continuous one-dimensional Schr\"odinger equation \cite{Bo:1946, Fa:1964, DeTr:1979, Si:1999,Ma:2011,Kravchenko1,Kravchenko2}, and the eigenvalue equation for a finite Jacobi matrix, or a semi-infinite Jacobi operator \cite{hochstadt1974construction,de1978numerically,GeSi:1997,Gi:TAMS2002,korotyaev2023inverse,MikhaylovMikhaylov,Mikhailov2025}. As the present article confirms, equation (\ref{discrete-2}) has features of both the finite and continuous structures; indeed our theoretical techniques draw on our previous work in both contexts \cite{Gi:TAMS2002,garnier2023effective, Gi:SIAMinpress,Gi:M3AS}.  Results in the present article complement previous works concerning periodic potentials \cite{Moerbeke1976,KoKu:2006,MiZh:2021}, multidimensional discrete problems \cite{isozaki2012inverse,isozaki2015inverse,Novikov1,Novikov2,Novikov3}, as well as the recent probabilistic analysis in the setting of one dimensional case \cite{garnier2023effective} and quasi-one dimensional case \cite{garnier2023surf,garnier2023robustness}.

Applicability of inverse theory to physical systems in an experimental setting requires computational results rather than merely theoretical. For instance, there is a considerable body of work concerned with periodic Jacobi operators, of which $\defJz$ is the most basic example. But the established proof that an $\ell_2$ perturbation of a periodic Jacobi operator is determined by scattering data, as presented in \cite[\S10.4]{teschl2000jacobi} for instance, does not translate easily into an efficient computational method. Similar remarks apply to \cite[Thm.~3.1]{volberg2002inverse}. An underlying motivation for the present article was to develop a simpler, more direct theoretical approach to the determination of $\oJ$ from scattering data that better lends itself to effective computation. In addition, new theoretical methods are required to decide whether scattering data for (\ref{discrete-2}) determines the perturbation length ${\nuL}$, since, as mentioned earlier, this question has not previously been addressed.

{\bf Outline:}
\S\ref{sec-main} details our main results from the theoretical perspective. Our main theoretical results therein are structured as follows. Theorem~\ref{thm-S-formula} expresses the scattering matrix for (\ref{discrete-2}) explicitly in terms of classical orthogonal polynomials. Theorem~\ref{thm-poles} analyzes precisely the zeros and poles of a certain linear fractional transformation of the reflection coefficient for (\ref{discrete-2})---which turns out to be crucial for determination of the perturbation length. Lastly, Theorem~\ref{thm-analytic-inversion} details a finite procedure to compute the unknown perturbed Schr\"odinger operator. The theoretical results are then translated to a practical algorithm in \S\ref{sec-computational}, whose implementation is illustrated via several numerical examples. The article concludes with a brief discussion in \S\ref{sec-discussion} of the key results and directions for further research based on certain questions left open in this article.  

\section{Main result\label{sec-main}}

\subsection{Matrix formulation and the $\mu$-function\label{sec-matrix}}

With
\begin{equation}\label{discrete-Un}
U_n(\sigma)\defeq\binom{u(n-1,\sigma)}{u(n,\sigma)},
\end{equation}
equation (\ref{discrete-2}) may be expressed in the form
\begin{equation}\label{discrete-3}
U_n=M_nU_{n+1}\quad\mbox{ where }\quad M_n\defeq\begin{pmatrix}\dfrac{a_n-\sigma^2}{b_{n-1}}&-\dfrac{b_n}{b_{n-1}}\\ 1&0\end{pmatrix}\qquad (n\in\integer).
\end{equation}
In terms of this notation, 
\begin{equation}\label{H}
U_1=HU_{{\nuL}+1}\quad\mbox{ where }\quad H\defeq\begin{pmatrix}h_{11}&h_{12}\\ h_{21}&h_{22}\end{pmatrix}=M_1\cdots M_{\nuL}. 
\end{equation}
Note by (\ref{discrete-3}) that entries $h_{ij}$ of $H$ are polynomial functions of $\sigma^2$, which according to the dispersion relation \eqref{dispersion}, has the form 
\begin{equation}\label{dispersion-2}
\sigma^2=2-z-\frac{1}{z}\quad\mbox{ where }\quad z=e^{ik}.
\end{equation}
The $\mu$-function for (\ref{discrete-2}) is defined in terms of solutions $u$ that are outgoing to the right, such as $f_1$ in (\ref{asymptotic}).  Fix notation as follows. 
Let $\mathbf{p}$ denote the projective mapping
\begin{equation}\label{projection}
\mathbf{p}:\complex^2\setminus\{0\}\rightarrow\complex\cup\{\infty\},\qquad\mathbf{p}\binom{z}{w}=\frac{z}{w}.
\end{equation}
Observe by (\ref{asymptotic}) and (\ref{H}) that for $u=f_1$,
\begin{equation}\label{pH}
\mathbf{p}U_1=\frac{u(0,\sigma)}{u(1,\sigma)}=\frac{1+R_2}{z+\frac{1}{z}R_2}=\mathbf{p}\left(HU_{{\nuL}+1}\right)=\frac{h_{11}+zh_{12}}{h_{21}+zh_{22}}\quad\mbox{ where }\quad z=e^{ik}.
\end{equation}
Taking reciprocals in (\ref{pH}) and expressing the $h_{ij}$ as functions of $\sigma^2=2-z-\frac{1}{z}$, yields
\begin{equation}\label{reciprocals}
\frac{u(1,\sigma)}{u(0,\sigma)}=\frac{z+\frac{1}{z}R_2}{1+R_2}=\frac{h_{21}(\sigma^2)+zh_{22}(\sigma^2)}{h_{11}(\sigma^2)+zh_{12}(\sigma^2)}\big|_{\sigma^2=\bigl(2-z-\frac{1}{z}\bigr)}\quad\mbox{ where }\quad z=e^{ik}.
\end{equation}
\begin{defn}\label{defn-mu}
The $\mu$-function for a perturbed Schr\"odinger operator (\ref{J}) is defined to be the rational function $\mu_{\oJ}:\complex\rightarrow\complex\cup\{\infty\}$ whose restriction to the circle $z=e^{ik}$ is given by the formula
\begin{equation}\label{mu-R}
\mu_{\oJ}\bigl(e^{ik}\bigr)=\frac{e^{ik}+e^{-ik}R_2(k)}{1+R_2(k)}. 
\end{equation}
\end{defn}
That $\mu_{\oJ}$ is a well-defined rational function of $z$ is evident from the right-hand expression in (\ref{reciprocals}).  Equation (\ref{reciprocals}) also shows that, given an outgoing-to-the-right nontrivial solution $u$ to (\ref{discrete-2}), the $\mu$-function for ${\oJ}$ is simply the extension to the complex plane of the ratio $u(1,\sigma)/u(0,\sigma)$.  In the unperturbed case ${\oJ}={\defJz}$, we have $f_1(n,\sigma)=e^{ikn}$ for all $n\in\integer$, corresponding to $R_2=0$ and $T_2=1$, and hence $\mu_{{\defJz}}(z)=z$.  

\begin{remark}
    
The definition of the $\mu$-function is analogous to that of the Titchmarsh-Weyl $m$-function in the continuous setting \cite{Si:1999},
\[
m(\sigma^2)=\frac{y^\prime(0,\sigma)}{y(0,\sigma)},
\]
where $y(x,\sigma)$ is a solution to the continuous one-dimensional Schr\"odinger equation that is outgoing to the right. Gesztesy and Simon in \cite[Prop.~2.3]{GeSi:1997} define an $m$-function for finite and semi-infinite Jacobi matrices to be what is effectively the negative of the left-most expression in (\ref{reciprocals}).
\end{remark}

\subsection{Orthogonal polynomials\label{sec-OPs}}

Referring to (\ref{J}), set 
\begin{equation}\label{J-n}
J_n\defeq
\begin{pmatrix}a_1&-b_1&0&\cdots&0\\
-b_1&a_2&-b_2&\ddots&\vdots\\
0&-b_2&\ddots&\ddots&0\\
\vdots&\ddots&\ddots&a_{n-1}&-b_{n-1}\\
0&\cdots&0&-b_{n-1}&a_n
\end{pmatrix},
\widetilde{J}_n\defeq
\begin{pmatrix}a_2&-b_2&0&\cdots&0\\
-b_2&a_3&-b_3&\ddots&\vdots\\
0&-b_2&\ddots&\ddots&0\\
\vdots&\ddots&\ddots&a_{n-1}&-b_{n-1}\\
0&\cdots&0&-b_{n-1}&a_n
\end{pmatrix},
\end{equation}
for $n\geq 1$ and $n\geq 2$, respectively.
Let $Q_0=S_0=1$, and, for $n\geq 1$, set 
\begin{equation}\label{QS}
Q_n(\lambda)\defeq\det(J_n-\lambda I)\quad\mbox{ and }\quad S_{n-1}(\lambda)\defeq\det(\widetilde{J}_n-\lambda I).
\end{equation}
Then $Q_n$ and $S_n$ are classical orthogonal polynomials on the line related to the corresponding sequences of finite Jacobi matrices.  
Indeed, writing $J_n=\Omega\Lambda \Omega^\top$ where $\Omega=(\omega_{ij})$ is orthogonal and $\Lambda={\rm diag}(\lambda_1,\ldots,\lambda_n)$, the sequence $Q_0,\ldots Q_{n-1}$ is orthogonal with respect to the pure point measure
\[
d\mu_n(x):=\sum_{j=1}^n\omega_{1j}^2\delta(x-\lambda_j). 
\]
Moreover, as $n\rightarrow\infty$ the measures $d\mu_n$ converge weakly to a probability measure $d\mu$ on $\mathbb{R}$ with respect to which the full sequence $Q_n$ $(n\geq 0)$ is orthogonal.
A similar statement holds for the sequence $S_0,\ldots S_{n-2}$; see \cite{de1978numerically}.
Define rescalings $q_n$ and $s_n$ by the formulas
\begin{equation}\label{rescaling}
q_n=\frac{Q_n}{b_1\cdots b_n}\quad\mbox{ and }\quad s_{n-1}=\frac{S_{n-1}}{b_1\cdots b_n},
\end{equation}
and set $s_{-1}=0$.  It follows from these definitions that, for $n\geq 2$, the recurrence 
\begin{equation}\label{qs-recurrence}
-b_{n-1}u(n-1,\sigma)+a_nu(n,\sigma)-b_nu(n+1,\sigma)=\sigma^2u(n,\sigma)
\end{equation}
is satisfied both with 
\begin{equation}\label{substitutions}
u(n,\sigma)=q_{n-1}(\sigma^2)\quad\mbox{ and }\quad u(n,\sigma)=s_{n-2}(\sigma^2)\qquad (n\geq 1). 
\end{equation}
A consequence in terms of the matrices $M_n$ defined in (\ref{discrete-3}) is that 
\begin{equation}\label{M-product}
M_1\cdots M_n=b_n\begin{pmatrix}q_n&-q_{n-1}\\ s_{n-1}&-s_{n-2}\end{pmatrix}\qquad(n\geq 1).
\end{equation}
Using the formula (\ref{discrete-3}) to take determinants in (\ref{M-product}), and noting $b_0=1$, yields
\begin{equation}\label{determinants}
q_{n-1}s_{n-1}-q_ns_{n-2}=1/b_n. 
\end{equation}

\subsection{Formulas in terms of orthogonal polynomials\label{sec-formulas}}

According to (\ref{J}) and (\ref{conditions}), $b_{\nuL}=1$ for a perturbed Schr\"odinger operator ${\oJ}$ having {\prtblen} ${\nuL}$. Therefore (\ref{H}) combined with (\ref{M-product}) yields 
\begin{equation}\label{H-formula}
H=\begin{pmatrix}q_{\nuL}&-q_{{\nuL}-1}\\ s_{{\nuL}-1}&-s_{{\nuL}-2}\end{pmatrix}.
\end{equation} 
Note by (\ref{determinants}) that 
\begin{equation}\label{H-determinant}
\det H=q_{{\nuL}-1}s_{{\nuL}-1}-q_{{\nuL}}s_{{\nuL}-2}=1.
\end{equation}
Given (\ref{H-formula}), application of equation (\ref{H}) to the special solutions $u=f_1$ and $u=f_2$ allows one to solve for the entries of the $S$-matrix, and hence also the $\mu$-function, in terms of orthogonal polynomials.  The resulting formulas are as follows. 
\begin{thm}\label{thm-S-formula}
Let ${\oJ}$ be a perturbed Schr\"odinger operator having {\prtblen} ${\nuL}$, and let polynomials $q_n, s_{n-1}$ $(n\geq0)$ be defined as in \S\ref{sec-OPs}.  Writing $z=e^{ik}$ and $\sigma^2=2-z-\frac{1}{z}$, the $S$-matrix for ${\oJ}$ as defined in (\ref{S}) has entries
\begin{equation}\label{S-formula}
\begin{split}
T_1(k)&=T_2(k)=\frac{1-z^2}{z^{\nuL}\bigl(q_{\nuL}(\sigma^2)-zq_{{\nuL}-1}(\sigma^2)-z(s_{{\nuL}-1}(\sigma^2)-zs_{{\nuL}-2}(\sigma^2))\bigr)}\\
R_2(k)&=z\frac{s_{{\nuL}-1}(\sigma^2)-zs_{{\nuL}-2}(\sigma^2)-z(q_{{\nuL}}(\sigma^2)-zq_{{\nuL}-1}(\sigma^2))}{q_{{\nuL}}(\sigma^2)-zq_{{\nuL}-1}(\sigma^2)-z(s_{{\nuL}-1}(\sigma^2)-zs_{{\nuL}-2}(\sigma^2))}\\
R_1(k)&=\frac{s_{{\nuL}-1}(\sigma^2)-\frac{1}{z}s_{{\nuL}-2}(\sigma^2)-\frac{1}{z}(q_{\nuL}(\sigma^2)-\frac{1}{z}q_{{\nuL}-1}(\sigma^2))}{z^{2{\nuL}-1}\bigl(q_{{\nuL}}(\sigma^2)-zq_{{\nuL}-1}(\sigma^2)-z(s_{{\nuL}-1}(\sigma^2)-zs_{{\nuL}-2}(\sigma^2))\bigr)}.\end{split}
\end{equation}
The $\mu$-function for ${\oJ}$ is given by the formula
\begin{equation}\label{mu-2}
\mu_{\oJ}(z)=\frac{s_{{\nuL}-1}(\sigma^2)-zs_{{\nuL}-2}(\sigma^2)}{q_{{\nuL}}(\sigma^2)-zq_{{\nuL}-1}(\sigma^2)}.
\end{equation}
\end{thm}
In light of  (\ref{H-determinant}), a direct evaluation, using formulas (\ref{S-formula}), confirms the expected identities corresponding to the fact that $S$ \eqref{S} is unitary: 
\begin{cor}\label{cor-identities}
$|T_1|^2+|R_1|^2=|T_2|^2+|R_2|^2=1, 
T_1\overline{R}_2+R_1\overline{T}_2=0. 
$
\end{cor}

\subsection{Analysis of $\mu_{\oJ}$\label{sec-analysis}}

Since $b_{\nuL}=1$ if ${\nuL}$ is the {\prtblen} for a perturbed Schr\"odinger operator ${\oJ}$, it follows from (\ref{rescaling}) that the formula (\ref{mu-2}) for $\mu_{\oJ}$ may be expressed in terms of the polynomials $Q_n,S_n$ defined in (\ref{QS}).  In detailed form,
\begin{equation}\label{mu-3}
\mu_{\oJ}(z)=\frac{b_1\cdots b_{\nuL}\bigl(s_{{\nuL}-1}(\sigma^2)-zs_{{\nuL}-2}(\sigma^2)\bigr)}{b_1\cdots b_{\nuL}\bigl(q_{{\nuL}}(\sigma^2)-zq_{{\nuL}-1}(\sigma^2)\bigr)}=
\frac{S_{{\nuL}-1}(\sigma^2)-zS_{{\nuL}-2}(\sigma^2)}{Q_{{\nuL}}(\sigma^2)-zQ_{{\nuL}-1}(\sigma^2)},
\end{equation}
where, as before, the extended dispersion relation is in force, whereby 
\begin{equation}\label{dispersion-3}
\sigma^2=2-z-\frac{1}{z}\qquad(z\in\complex). 
\end{equation}
By the same token, it follows from (\ref{H-determinant}) that 
\begin{equation}\label{H-determinant-2}
Q_{{\nuL}-1}S_{{\nuL}-1}-Q_{{\nuL}}S_{{\nuL}-2}=b_1\cdots b_{\nuL}=b_1\cdots b_{{\nuL}-1}>0.
\end{equation}
In addition, the definition of the {\prtblen} implies that $a_{\nuL}\neq2$ or $b_{{\nuL}-1}\neq 1$.  
These observations facilitate precise determination of the numbers of zeros and poles of $\mu_{\oJ}$, as follows.  
\begin{thm}\label{thm-poles}
Let perturbed Schr\"odinger operator ${\oJ}$ have {\prtblen} ${\nuL}\geq1$.  Then $\mu_{\oJ}(z)$ has a simple zero at $z=0$.  If $a_{\nuL}\neq 2$ then $\mu_{\oJ}$ has $2{\nuL}-2$ zeros and $2{\nuL}-1$ poles. Otherwise $a_{\nuL}=2$ and $b_{{\nuL}-1}\neq 1$, in which case ${\nuL}\geq 2$ and $\mu_{\oJ}$ has $2{\nuL}-3$ zeros and $2{\nuL}-2$ poles. 
\end{thm}
\begin{proof}
Referring to (\ref{mu-3}), set 
\begin{equation}\label{AB}
A(z)=S_{{\nuL}-1}(\sigma^2)-zS_{{\nuL}-2}(\sigma^2)\quad\mbox{ and }\quad B(z)=Q_{{\nuL}}(\sigma^2)-zQ_{{\nuL}-1}(\sigma^2). 
\end{equation}
Consider first $B$. By (\ref{QS}), 
\[
Q_{\nuL}(\sigma^2)=(a_{\nuL}-\sigma^2)Q_{{\nuL}-1}(\sigma^2)-b_{{\nuL}-1}^2Q_{{\nuL}-2}(\sigma^2),
\]
and so 
\begin{equation}\label{B-1}
\begin{split}
B(z)&=(a_{\nuL}-\sigma^2-z)Q_{{\nuL}-1}(\sigma^2)-b_{{\nuL}-1}^2Q_{{\nuL}-2}(\sigma^2)\\
&=(a_{\nuL}-2+\textstyle\frac{1}{z})Q_{{\nuL}-1}(2-z-\frac{1}{z})-b_{{\nuL}-1}^2Q_{{\nuL}-2}(2-z-\frac{1}{z}).
\end{split}
\end{equation}
In general 
\begin{equation}\label{Q-general}
Q_n(\sigma^2)=(-\sigma^2)^n+O(\sigma^{2(n-1)})=(z+\textstyle\frac{1}{z}-2)^n+O(z^{n-1}).  
\end{equation}
Therefore the largest negative power of $z$ in the expansion of $B(z)$ in (\ref{B-1}) is $z^{-{\nuL}}$, coming from the product of $1/z$ with the leading term  $(z+\frac{1}{z}-2)^{{\nuL}-1}$ of $Q_{{\nuL}-1}$.  Since this is the unique occurrence of $z^{-{\nuL}}$, it follows that $z^{\nuL} B(z)$ is a polynomial whose constant term is 1.  The degree of the polynomial $z^{\nuL} B(z)$ depends on $a_{\nuL}$.  If $a_{\nuL}-2\neq 0$, then it follows from (\ref{B-1}) and (\ref{Q-general}) that $B$ has the form 
\[
B(z)=(a_{\nuL}-2)(z+\textstyle\frac{1}{z}-2)^{{\nuL}-1}+O(z^{{\nuL}-2}),
\]
whereby the degree of $z^{\nuL} B(z)$ is $2{\nuL}-1$.  On the other hand, if $a_{\nuL}-2=0$, then it follows from (\ref{B-1}) and (\ref{Q-general}) that $B(z)$ has the form
\[
B(z)=\textstyle\frac{1}{z}(z+\frac{1}{z}-2)^{{\nuL}-1}-b_{{\nuL}-1}^2(z+\frac{1}{z}-2)^{{\nuL}-2}+O(z^{{\nuL}-3})=(1-b_{{\nuL}-1}^2)z^{{\nuL}-2}+O(z^{{\nuL}-3}),
\]
in which case $z^{\nuL} B(z)$ has precise degree $2{\nuL}-2$.  

Exactly the same reasoning applies to $A(z)$ and the orthogonal polynomials $S_n(\sigma^2)$, to yield that $z^{{\nuL}-1}A(z)$ is a polynomial of constant term 1 and degree $2{\nuL}-3$ if $a_{\nuL}-2\neq0$, or degree $2{\nuL}-4$ if $a_{\nuL}-2=0$.  Therefore 
\begin{equation}\label{mu-4}
\mu_{\oJ}(z)=\frac{A(z)}{B(z)}=z\frac{z^{{\nuL}-1}A(z)}{z^{\nuL} B(z)}
\end{equation}
has a simple zero at $z=0$.  It remains to show that the polynomials $z^{{\nuL}}A(z)$ and $z^{{\nuL}}B(z)$ have no common zero.  Suppose on the contrary that $z_0$ is a common zero of $z^{{\nuL}}A(z)$ and $z^{{\nuL}}B(z)$.  Then $z_0\neq 0$, since 
\[
\left.z^{\nuL} B(z)\right|_{z=0}=1. 
\]
It follows therefore that $A(z_0)=B(z_0)=0$, and hence by (\ref{AB}) that 
\[
S_{{\nuL}-1}(2-z_0-\textstyle\frac{1}{z_0})=z_0S_{{\nuL}-2}(2-z_0-\frac{1}{z_0})\quad\mbox{ and }\quad Q_{{\nuL}}(2-z_0-\textstyle\frac{1}{z_0})=z_0Q_{{\nuL}-1}(2-z_0-\frac{1}{z_0}).
\]
But then 
\[
0=z_0\bigl(Q_{{\nuL}-1}(2-z_0-\frac{1}{z_0})S_{{\nuL}-1}(2-z_0-\textstyle\frac{1}{z_0})-Q_{{\nuL}}(2-z_0-\textstyle\frac{1}{z_0})S_{{\nuL}-2}(2-z_0-\frac{1}{z_0}),
\]
contradicting (\ref{H-determinant-2}).  This proves polynomials $z^{{\nuL}}A(z)$ and $z^{{\nuL}}B(z)$ have no common zero. It follows by (\ref{mu-4}) that the zeros of $\mu_{\oJ}(z)$ are 0 together with the zeros of the polynomial $z^{{\nuL}-1}A(z)$, and the poles of $\mu_{\oJ}(z)$ are the zeros of $z^{{\nuL}}B(z)$, the precise number depending on whether or not $a_{\nuL}=2$, as detailed above. 
\end{proof}

\begin{cor}\label{cor-poles}
The {\prtblen} ${\nuL}$ for a perturbed Schr\"odinger operator ${\oJ}$ is given by the formula 
$
{\nuL}=\lceil (p+1)/2\rceil,
$
where $p$ denotes the number of poles of $\mu_{\oJ}$. 
\end{cor}

The next objective is to express $\mu_{\oJ}$ as a composition of linear fractional transformations from which the perturbed values $a_n$ and $b_n$ can be extracted.  Let
\begin{equation}\label{defmatD}
D\defeq\begin{pmatrix}1&0\\ 0&-1\end{pmatrix}. 
\end{equation}
It follows from (\ref{reciprocals}) that 
\begin{equation}\label{mu-5}
\mu_{\oJ}(z)=\mathbf{p}\left(G\binom{z}{1}\right)\quad\mbox{ where }\quad G\defeq\begin{pmatrix}h_{22}&h_{21}\\ h_{12}& h_{11}\end{pmatrix}=D(H^{-1})^\top D,
\end{equation}
using the definition of $\mathbf{p}$ \eqref{projection}.
Since $H=M_1\cdots M_{\nuL}$, the matrix $G$ in (\ref{mu-5}) factors as
\begin{equation}\label{G-factors}
G=G_1\cdots G_{\nuL}\quad\mbox{ where }\quad G_n=D(M_n^{-1})^\top D=\begin{pmatrix}0&\dfrac{b_{n-1}}{b_n}\\ -1&\dfrac{a_n-\sigma^2}{b_n}\end{pmatrix}\qquad(1\leq n\leq {\nuL}),
\end{equation}
where $\sigma^2=2-z-\frac{1}{z}$, as usual.  Define associated linear fractional transformations
\begin{equation}\label{g}
g_{n,z}(\xi)=\mathbf{p}\left(G_n\binom{\xi}{1}\right)=\frac{b_{n-1}z}{z(a_n-2)+z^2+1-b_nz\xi}\qquad(1\leq n\leq {\nuL};\thickspace z,\xi\in\complex).
\end{equation}
Then (\ref{mu-5}) may be expressed in the form
\begin{equation}\label{mu-6}
\mu_{\oJ}(z)=g_{1,z}\circ\cdots\circ g_{{\nuL},z}(z). 
\end{equation}
The foregoing definitions lead directly to the following recursion. 
\begin{lem}\label{lem-recurrence}
Let perturbed Schr\"odinger operator ${\oJ}$ of the form (\ref{J}) have {\prtblen} ${\nuL}\geq 1$. Set $\varphi_{{\nuL}+1}(z)=z$, and conisder the definition
\[
\varphi_n(z)\defeq b_{n-1}g_{n,z}\circ\cdots\circ g_{{\nuL},z}(z)\qquad(1\leq n\leq {\nuL}),
\]
where $g_{n,z}$ is defined as in (\ref{g}).  Then $\varphi_1=\mu_{\oJ}$, and 
\begin{equation}\label{recurrence}
\varphi_n(z)=\frac{b_{n-1}^2z}{z(a_n-2)+z^2+1-z\varphi_{n+1}(z)}\qquad(1\leq n\leq {\nuL}). 
\end{equation}
\end{lem}

\subsection{Inverse scattering\label{sec-inverse}}
Corollary~\ref{cor-poles} and Lemma~\ref{lem-recurrence} justify the following theoretical inversion procedure, formulated in terms of the reflection coefficient $R_2$.  The procedure centres on formula (\ref{recurrence}), from which $b_{n-1}$, $a_n$ and $\varphi_{n+1}$ can be extracted one after another.

\begin{thm}\label{thm-analytic-inversion}
Let \[S=\begin{pmatrix}T_1&R_2\\ R_1&T_2\end{pmatrix}\] be the scattering matrix for a perturbed Schr\"odinger operator ${\oJ}$ of the form (\ref{J}) that is a finitely supported perturbation of ${\defJz}$.  Given $R_2(k)$ for all $k\in(-\pi,\pi]$, the following procedure determines ${\oJ}$.
\begin{enumerate}
\item If $R_2=0$, then ${\oJ}={\defJz}$ and we are done.  Otherwise, let $\varphi_1(z)$ denote the unique rational extension to $z\in\complex$ of the formula
\[
\varphi_1(z)=\frac{z+\frac{1}{z}R_2(k)}{1+R_2(k)}\qquad(z=e^{ik};\thickspace -\pi<k\leq\pi).
\]
Let $p$ denote the number of poles of $\varphi_1$.  The {\prtblen} for ${\oJ}$ is \[{\nuL}=\lceil(p+1)/2\rceil.  \]
\item For $n=1$ to ${\nuL}$, given $\varphi_n$:
\begin{enumerate}[label=2.\arabic*.]
\item determine $b_{n-1}>0$ by the formula
$b_{n-1}^2=\left.\frac{1}{z}\varphi_n(z)\right|_{z=0};$
\item given $b_{n-1}^2$, set
$a_n=\left.\frac{b_{n-1}^2}{\varphi_n(z)}+2-z-\frac{1}{z}\right|_{z=0};$
\item given $b_{n-1}^2$ and $a_n$, set 
$\varphi_{n+1}(z)=a_n-2+z+\frac{1}{z}-\frac{b_{n-1}^2}{\varphi_n(z)}.$\\
End.
\end{enumerate}
\item Output $b_1,\ldots,b_{{\nuL}-1}$ and $a_1,\ldots,a_{\nuL}$. 
\end{enumerate}
\end{thm}
Thus we have a finite procedure to determine ${\oJ}$ from its scattering data.  The procedure entails counting the number of poles of a rational function given its restriction to the unit circle, and, in steps 2.1 and 2.2, evaluation of various rational functions defined in terms of data given on the unit circle at the point $z=0$.  Translation of Theorem~\ref{thm-analytic-inversion} into a computationally effective inversion scheme requires that these steps be carried out numerically.  This involves some analysis of poles of rational functions, carried out in the next section.  

Note the following immediate consequence of Theorem~\ref{thm-analytic-inversion} concerning reflectionless perturbed Schr\"odinger operators, i.e. perturbed Schr\"odinger operators for which $R_2=0$; see \cite{jakvsic2014note}. 
\begin{cor}\label{cor-reflectionless}
A perturbed Schr\"odinger operator $\oJ$ satisfying (\ref{J}) and (\ref{conditions}) is reflectionless only in the trivial case $\oJ=\defJz$. 
\end{cor}

In the present section, the notation $\complex^d$, for example, indicates column vectors, and is synonymous with $\complex^{d\times 1}$. It is convenient to fix the form of a grid of $d$ sample points in $(-\pi,\pi)$ as
\begin{equation}\label{sample-grid}
k_j=-\pi-\frac{\pi}{d}+j\frac{2\pi}{d}\qquad(1\leq j\leq d).
\end{equation}
Corresponding points on $S^1$ are
\begin{equation}\label{sample-grid-2}
z_j=e^{i\left(-\pi-\frac{\pi}{d}+j\frac{2\pi}{d}\right)}=-e^{-i\pi/d}(e^{2\pi i/d})^j\qquad(1\leq j\leq d).
\end{equation}
It will be assumed that scattering data for a perturbed Schr\"odinger operator $\oJ$ is given in the form of the reflection coefficient $R_2(k_j)$, where $k_j$ ranges over an evenly-spaced grid of sample points of the form (\ref{sample-grid}).  In addition, it will be assumed that the order of magnitude of the perturbation length $L$ of $\oJ$ is known to be less than a given $N$, where $N$ is in turn less than the grid size $d$ ($N$ is assumed to be known based on physical considerations, say, while $d$ is inherent in the data itself).  Moreover, it will be assumed that values $L,N,d$ are different orders of magnitude.  In practice we require
\begin{equation}\label{relative-sizes}
2L-1\leq N\quad\mbox{ and }\quad 2N-1\leq d.
\end{equation}

\section{Numerical computations}
\label{sec-computational}

Translation of Theorem~\ref{thm-analytic-inversion} into a practical algorithm requires some additional analysis of poles of rational functions, and corresponding numerical algorithms to count poles (Algorithm~\ref{alg-pole-count} {\sc PoleCount}) and correct for the associated residues (Algorithn-\ref{alg-pole-correction} {\sc PoleCorrection}). The analysis of poles and related algorithms are detailed in Appendices~\ref{sec-poles} and \ref{sec-algorithms}. With these in hand, Theorem~\ref{thm-analytic-inversion} gives rise to the following algorithm. Here $v$ denotes a vector of $d$ equally-spaced samples on the unit circle of $R_2$; $N$ is a preset upper bound on the perturbation length; and $\epsilon$ denotes a prescribed tolerance (applicable to a singular value count). 

\floatname{algorithm}{Main Algorithm}
\renewcommand{\thealgorithm}{INVR} 
\begin{algorithm}[htb!]
\caption{Inversion of $R_2$ sampled on $(-\pi,\pi)$}
\label{alg-reflection-inversion}
\begin{algorithmic}[1]
\Procedure{ReflInv}{$v,N,\epsilon$}
\State {\bf Input} $(v,N,\epsilon)\in\complex^d\times\{1,\ldots,\lfloor (d+1)/2\rfloor\}\times\real_+$
\State $z_m\gets -e^{-i\pi/d}( e^{2\pi i/d})^m\quad(1\leq m\leq d)$
\State $\mu\gets\left(\frac{z_1+\frac{1}{z_1}v_1}{1+v_1},\ldots,\frac{z_d+\frac{1}{z_d}v_d}{1+v_d}\right)^\top$
\State ${\nuLg}\gets\lceil\left(\mbox{\sc PoleCount}(v,N,\epsilon)_1+\mbox{\sc PoleCount}(\overline{v},N,\epsilon)_1+1\right)/2\rceil$ \label{overcount-line}
\State $\tau\gets\left(\frac{1}{\mu_1}+2-z_1-\frac{1}{z_1},\ldots,\frac{1}{\mu_d}+2-z_d-\frac{1}{z_d}\right)^\top$
\State $a_1\gets\mbox{\sc PoleCorrection}(\tau,N,\epsilon)+\frac{1}{d}\sum_{m=1}^d\tau_m$
\State $\phi_m\gets a_1-2+z_m+\frac{1}{z_m}-\frac{1}{\mu_m}\quad(1\leq m\leq d)$
\State $\tau_m\gets {\phi_m}/{z_m}\quad(1\leq j\leq d)$
\For{$j=1$ to ${\nuLg}-1$}
\State $b_j\gets\left(\mbox{\sc PoleCorrection}(\tau,N,\epsilon)+\frac{1}{d}\sum_{m=1}^d\tau_m\right)^{1/2}$
\State $\tau_m\gets \frac{b_j^2}{\phi_m}+2-z_m-\frac{1}{z_m}\quad(1\leq m\leq d)$
\State $a_{j+1}\gets\mbox{\sc PoleCorrection}(\tau,N,\epsilon)+\frac{1}{d}\sum_{m=1}^d\tau_m$
\State $\phi_m\gets a_{j+1}-2+z_m+\frac{1}{z_m}-\frac{b_j^2}{\phi_m}\quad(1\leq m\leq d)$
\State $\tau_m\gets \phi_m/z_m\quad(1\leq m\leq d)$
\EndFor
\State {\bf Output} $(a_1,\ldots,a_{\nuLg},b_1,\ldots,b_{{\nuLg}-1})\in\complex^{2{\nuLg}-1}$
\EndProcedure
\end{algorithmic}
\end{algorithm}
Algorithm~\ref{alg-reflection-inversion} is easily implemented in any commonly-used computing environment, since it reduces the inverse scattering problem to standard numerical methods: SVD, the generalized eigenvalue problem, and the least-squares solution to a linear system.

The following examples were produced by implementation of Algorithm~\ref{alg-reflection-inversion} in \emph{Mathematica} with $(N,\epsilon)=(100,10^{-10})$ in all cases. Perturbations of size $S$ are generated from a uniformly distributed random sample in the interval $[-S,S]$; these perturbations are added to $\defJz$ over the perturbation length. (Thus, $\max|\defJz-\oJ|\leq S$.) Recall that $N$ serves as an upper bound on the perturbation length, and $\epsilon$ is the tolerance for singular value count. The examples serve as a proof of concept of the theoretical results encapsulated in Corollary~\ref{cor-poles}, Theorem~\ref{thm-analytic-inversion} and Proposition~\ref{prop-poles}.

\begin{figure}[htb!]
\centering
\begin{subfigure}[b]{0.4\textwidth}
    \centering
    \includegraphics[width=\textwidth]{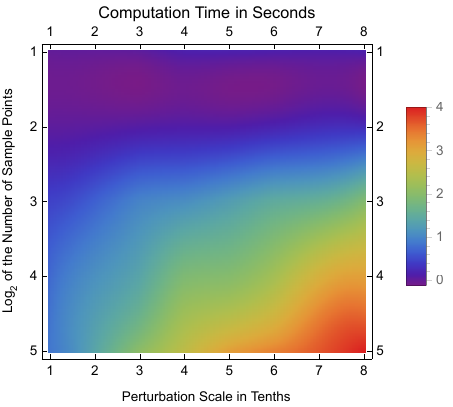}
    \caption{Timing for length 4 perturbations.}
    \label{fig-0a}
\end{subfigure}
\begin{subfigure}[b]{0.4\textwidth}
    \centering
    \includegraphics[width=\textwidth]{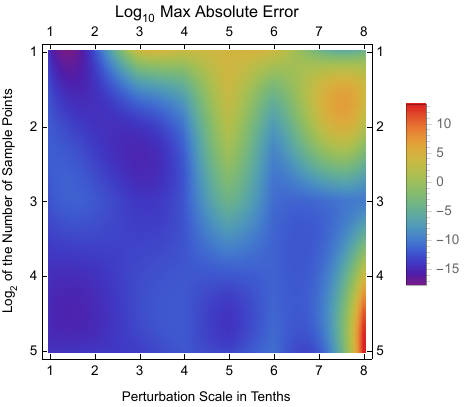}
    \caption{Maximum absolute error for length 4 perturbations.}
    \label{fig-0b}
\end{subfigure}
\vspace{0.5cm}
\begin{subfigure}[b]{0.4\textwidth}
    \centering
    \includegraphics[width=\textwidth]{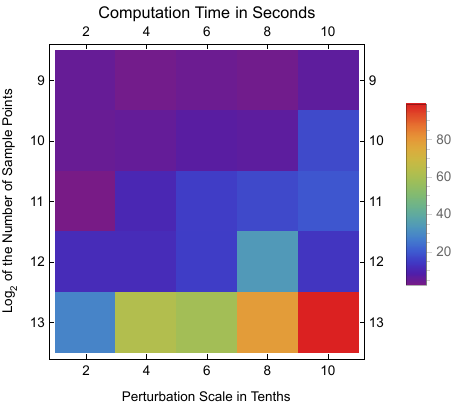}
    \caption{Timing for length 8 perturbations.}
    \label{fig-01a}
\end{subfigure}
\begin{subfigure}[b]{0.4\textwidth}
    \centering
    \includegraphics[width=\textwidth]{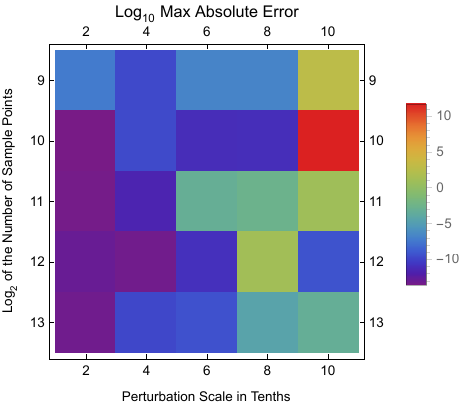}
    \caption{Maximum absolute error for length 8 perturbations.}
\label{fig-01b}
\end{subfigure}
\caption{
Timing and maximum absolute error as functions of perturbation size and number of sample points for perturbation lengths 4 (top row) and 8 (bottom row). The length-4 results are averages over 40 runs with randomly generated perturbations. 
}
\label{fig:timing_error_comparison}
\end{figure}

Computational time (for the implementation of Algorithm~\ref{alg-reflection-inversion}) and accuracy (maximum absolute error; refer, for example, third column of Table \ref{table-examples}) depend on perturbation length $L$, perturbation size $S$, and the number of sample points $d$ on the circle. Figures \ref{fig:timing_error_comparison}
illustrate this dependency for respective perturbation lengths four and eight.  In the examples that follow, the number of sample points is either $10^3$ or $10^4$, the larger number being necessary in the case of larger perturbation length, as reflected in Figure~\ref{fig-0a}-\ref{fig-0b}.

\begin{table}[htb!]
\caption{Exact vs.~computed perturbations for Examples 1--4.}
\label{table-examples}
\centering
\begin{tabular}{|p{.4cm}|p{3.55cm}|p{.4cm}|p{3.55cm}|l|}
\hline
Ex. & $C_{a,b}$ & $\nu L_g$ & $\mbox{\sc ReflInv}(R_2,N,\epsilon)$ & Absolute error \\
\hline\hline

1 &
$\begin{pmatrix}
1.89266 & 1. \\
1.80067 & 1. \\
1.98101 & 1. \\
2.01661 & 1. \\
2.04581 & 1. \\
1.89497 & 1. \\
1.94019 & 1. \\
2.14114 &
\end{pmatrix}$
&
8
&
$\begin{pmatrix}
1.89266 & 1. \\
1.80067 & 1. \\
1.98101 & 1. \\
2.01661 & 1. \\
2.04581 & 1. \\
1.89497 & 1. \\
1.94019 & 1. \\
2.14114 &
\end{pmatrix}$
&
$10^{-16}\times
\begin{pmatrix}
0      & 1.1102\\
13.323 & 2.2204\\
0      & 1.1102\\
4.4409 & 2.2204\\
8.8818 & 3.3307\\
4.4409 & 4.4409\\
13.323 & 6.6613\\
4.4409 & \cdots
\end{pmatrix}$
\\
\hline

2 &
$\begin{pmatrix}
2. & 0.975627\\
2. & 1.08951\\
2. & 1.07743\\
2. & 0.907199\\
2. & 0.969677\\
2. & 1.09471\\
2. & 0.958727\\
2. &
\end{pmatrix}$
&
8
&
$\begin{pmatrix}
2. & 0.975627\\
2. & 1.08951\\
2. & 1.07743\\
2. & 0.907199\\
2. & 0.969677\\
2. & 1.09471\\
2. & 0.958727\\
2. &
\end{pmatrix}$
&
$10^{-16}\times
\begin{pmatrix}
2.2204 & 0\\
0 & 6.6613\\
6.6613 & 6.6613\\
8.8818 & 2.2204\\
4.4409 & 8.8818\\
13.323 & 4.4409\\
13.323 & 59.952\\
0 & \cdots
\end{pmatrix}$
\\
\hline

3 &
$\begin{pmatrix}
2.17037 & 0.922992\\
2.16532 & 1.00227\\
2.03169 & 1.05907\\
2.02329 & 1.0594\\
2.11997 & 0.91648\\
2.11719 & 0.912865\\
1.97313 & 1.02128\\
2.042 & 1.02982\\
1.86884 & 1.0648\\
2.18595 & 1.07725\\
2.13613 & 1.0392\\
1.95432 & 0.985749\\
1.91324 & 1.05941\\
1.96332 & 1.06087\\
2.19366 & 0.975405\\
2.05174 &
\end{pmatrix}$
&
16
&
$\begin{pmatrix}
2.17037 & 0.922992\\
2.16532 & 1.00227\\
2.03169 & 1.05907\\
2.02329 & 1.0594\\
2.11997 & 0.91648\\
2.11719 & 0.912865\\
1.97313 & 1.02128\\
2.042 & 1.02982\\
1.86884 & 1.0648\\
2.18595 & 1.07725\\
2.13613 & 1.0392\\
1.95432 & 0.985749\\
1.91324 & 1.05941\\
1.96332 & 1.06087\\
2.19366 & 0.975405\\
2.05174 &
\end{pmatrix}$
&
$10^{-10}\times
\begin{pmatrix}
0.000022204 & 0.0025213\\
0.03527 & 0.041018\\
0.079607 & 0.10678\\
0.063203 & 0.16843\\
0.0032552 & 0.34264\\
0.17768 & 0.97977\\
0.41670 & 2.0208\\
0.0063549 & 3.0191\\
3.1583 & 3.5709\\
4.8290 & 2.1526\\
1.7588 & 0.38595\\
2.9765 & 0.89504\\
2.0257 & 1.5169\\
6.3780 & 0.14580\\
1.2429 & 2.1021\\
2.1841 & \cdots
\end{pmatrix}$
\\
\hline

4 &
$\begin{pmatrix}
0.209093 & 1.26857\\
2.64167 & 1.15019\\
1.67696 & 1.73568\\
0.0983645 & 0.588488\\
1.6473 & 0.82317\\
3.04075 & 0.909752\\
1.12573 & 1.99893\\
3.3695 &
\end{pmatrix}$
&
10
&
$\begin{pmatrix}
0.209093 & 1.26857\\
2.64167 & 1.15019\\
1.67696 & 1.73568\\
0.0983645 & 0.588488\\
1.6473 & 0.82317\\
3.04075 & 0.909752\\
1.12572 & 1.99893\\
3.3695 & 1.\\
2. & 1.\\
2. &
\end{pmatrix}$
&
$10^{-7}\times
\begin{pmatrix}
0.000013329 & 0.000065645\\
0.00025521 & 0.00091489\\
0.0037124 & 0.0041904\\
0.0035233 & 0.016664\\
0.19555 & 0.41296\\
1.8310 & 1.7834\\
15.108 & 10.587\\
7.4079 & 5.1270\\
3.5902 & 0.62871\\
0.42730 & \cdots
\end{pmatrix}$
\\
\hline

\end{tabular}
\end{table}

\subsection{Example 1: a Schr\"odinger potential\label{sec-examples1}}
For a Schr\"odinger operator (\ref{discrete}) the coefficients $b_n$ in (\ref{J}) are all 1; just the $a_n$ are perturbed. We consider the case of a potential $v$ whose support has length $\nuL=8$, with the reflection coefficient being sampled on $d=10^3$ values of $k$. The data $R_2(k_j)$ and corresponding $\mu$-function $\mu_{\oJ}(k_j)$ are plotted in Figure~\ref{fig:combined_R_mu}(a), (e). The first example in Table~\ref{table-examples} compares the output of Algorithm~\ref{alg-reflection-inversion} to the actual values; the value of $\ell$ from line~\ref{overcount-line} of Algorithm~\ref{alg-reflection-inversion} is indicated in the second column. In principle $\ell\geq \nuL$ may be an overcount, since $\mu_{\oJ}$ may have poles on the circle, but in this case $\ell=\nuL$. Note that reconstruction is accurate essentially to machine precision.

In Table \ref{table-examples}, the first column represents
$C_{a,b}$ defined by
\begin{equation}\label{table-col-1}
C_{a,b}\defeq\begin{pmatrix}a_1&b_1\\ \vdots&\vdots\\ &b_{{\nuLg}-1}\\ a_{\nuLg}&\end{pmatrix}.
\end{equation}

\subsection{Example 2: perturbation of just the $b_n$\label{sec-examples2}}
This example is a sort of dual to the previous one: instead of perturbing just the $a_n$, we perturb just the $b_n$. As in the previous case, $\nuL=8$ and $d=10^3$. The symmetry of $\mu_{\oJ}$ across both axes manifest in Figure~\ref{fig:combined_R_mu}(f) is characteristic of $a_n\equiv 2$. The second example in Table~\ref{table-examples} indicates, the reconstruction is again essentially perfect. For the definition of $C_{a,b}$, recall \eqref{table-col-1}.

\subsection{Example 3: a general perturbation\label{sec-examples3}}
In this example Algorithm~\ref{alg-reflection-inversion} is given $R_2$ sampled on $d=10^3$ points for a general perturbation of length $\nuL=16$.  The third example in Table~\ref{table-examples} indicates a slight loss of accuracy in the reconstruction with increasing index $n$, which is what one would expect given the recursive nature of Algorithm~\ref{alg-reflection-inversion}. 

\subsection{Example 4: poles of $\mu_{\oJ}$ on $S^1$\label{sec-examples4}}
\begin{figure*}[htb!]
\centering
\begin{tabular}{c}
\includegraphics[width=.215\textwidth]{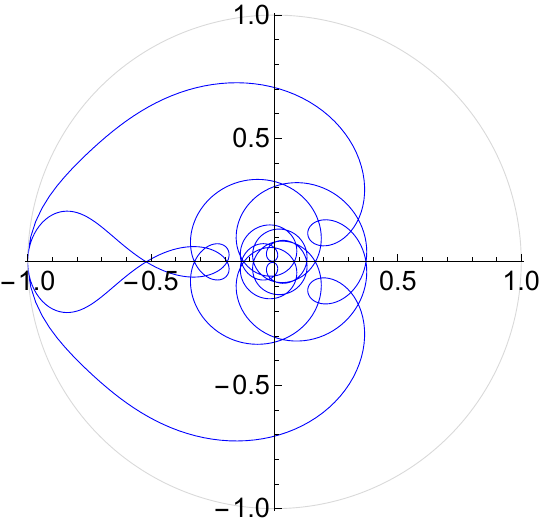}\\[-2mm]
(a)
\end{tabular}\hspace{-1mm}
\begin{tabular}{c}
\includegraphics[width=.215\textwidth]{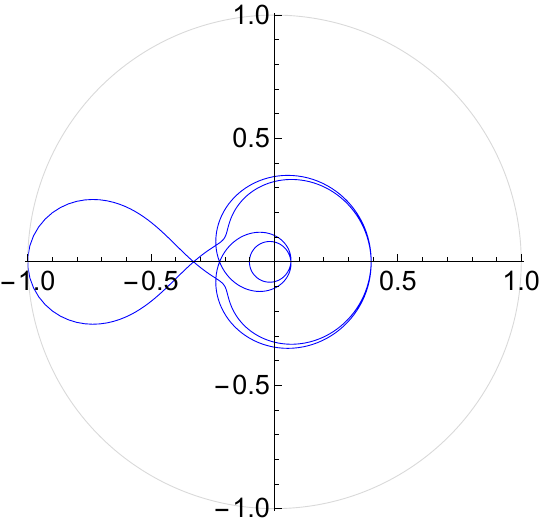}\\[-2mm]
(b)
\end{tabular}\hspace{-1mm}
\begin{tabular}{c}
\includegraphics[width=.215\textwidth]{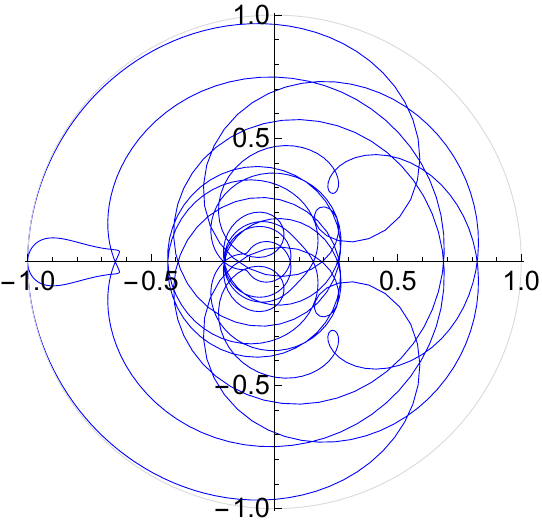}\\[-2mm]
(c)
\end{tabular}\hspace{-1mm}
\begin{tabular}{c}
\includegraphics[width=.215\textwidth]{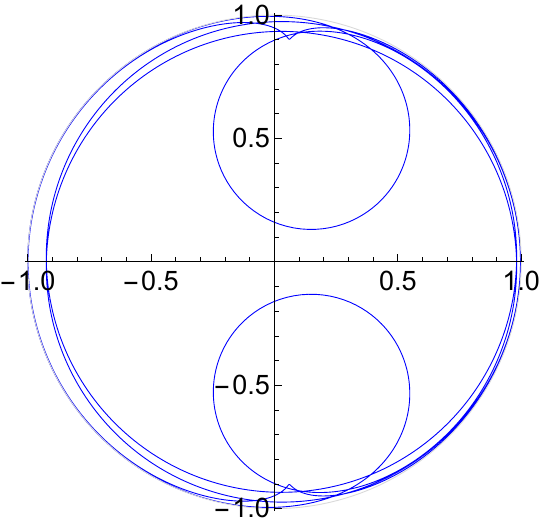}\\[-2mm]
(d)
\end{tabular}

\vspace{1mm}
\begin{tabular}{c}
\includegraphics[width=.215\textwidth]{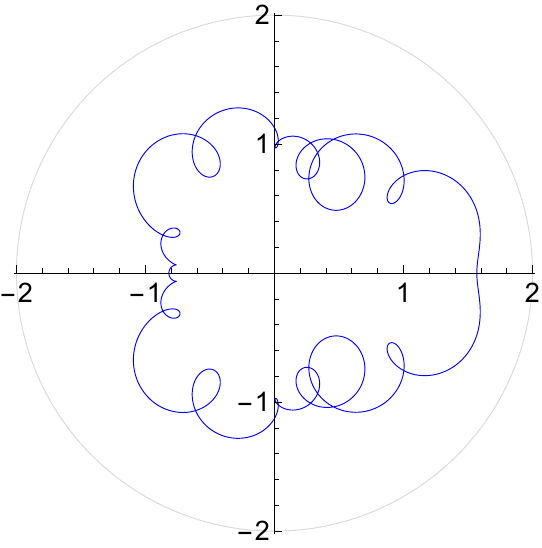}\\[-2mm]
(e)
\end{tabular}\hspace{-1mm}
\begin{tabular}{c}
\includegraphics[width=.215\textwidth]{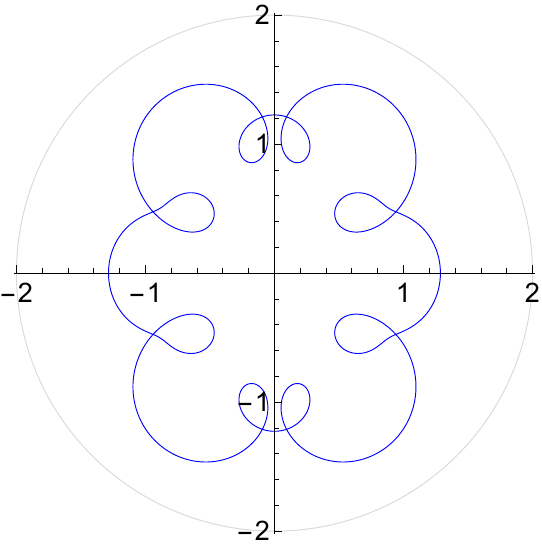}\\[-2mm]
(f)
\end{tabular}\hspace{-1mm}
\begin{tabular}{c}
\includegraphics[width=.215\textwidth]{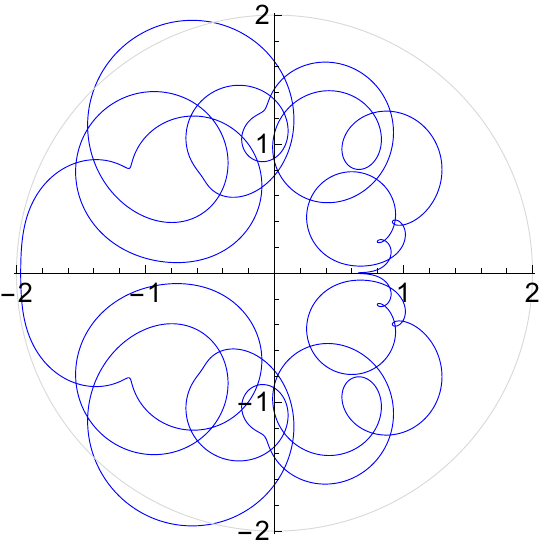}\\[-2mm]
(g)
\end{tabular}\hspace{-1mm}
\begin{tabular}{c}
\includegraphics[width=.215\textwidth]{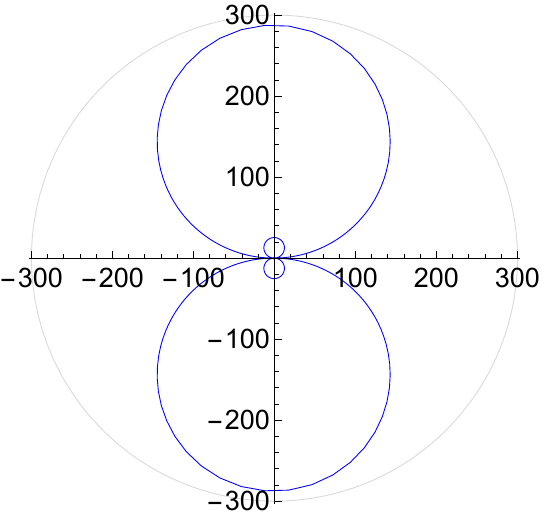}\\[-2mm]
(h)
\end{tabular}

\caption{
First row: $R_2(k_j)$ $(1\leq j\leq d)$ for perturbations 0.03 (length 8), 0.09 (length 8), 0.03 (length 16), and 0.20 (length 10), respectively.
Second row: Corresponding $\mu_J(k_j)$. 
}
\label{fig:combined_R_mu}
\end{figure*}

In this example, the reflection coefficient is sampled on $d=10^4$ points for a perturbed Schr\"odinger operator of perturbation length $\nuL=8$. The perturbation is such that $\mu_{\oJ}$ has poles on the unit circle (numerically, i.e., within machine precision), with the result that line~\ref{overcount-line} of Algorithm~\ref{alg-reflection-inversion} gives an overcount of $\ell=L+2$, and hence reconstructs $a_n$ and $b_n$ at two points beyond the actual perturbation. But these are correctly reconstructed as the unperturbed values, as the fourth example in Table~\ref{table-examples} shows, although there is a greater loss of accuracy than in the previous examples. The existence of poles on (or very near) $S^1$ is manifest in the comparatively large values of $\mu_{\oJ}$ appearing in Figure~\ref{fig:combined_R_mu}(h). In general, poles on the unit circle create an inherent instability in the reconstruction.  

\section{Discussion\label{sec-discussion}}

\subsection{Orthogonal polynomials and hybrid features}
A key technical ingredient of the present article is to isolate and analyze the polynomials
\begin{equation}\label{key-polynomials}
{S_{{\nuL}-1}(\sigma^2)-zS_{{\nuL}-2}(\sigma^2)}\quad\mbox{ and }\quad {Q_{{\nuL}}(\sigma^2)-zQ_{{\nuL}-1}(\sigma^2)}
\end{equation}
occurring in the formulas for $R_2$ and $\mu_{\oJ}$ in (\ref{S-formula}) and (\ref{mu-3}). 
These are polynomials in the variable $z=e^{ik}$, since $\sigma^2=2-z-1/z$, and are evaluated on the unit circle to yield the given data $R_2(k)$. But they are built from successive pairs $S_{\nuL-2}, S_{\nuL-1}$ and $Q_{\nuL-1},Q_{\nuL}$ of classical orthogonal polynomials on the line, as described in \S\ref{sec-OPs}. The central analytic object, which underpins both computation of $\nuL$ and determination of $\oJ$, is hybrid: a polynomial on the circle built from polynomials on the line. Whereas orthogonal polynomials on the line are central to the inverse theory of finite Jacobi operators \cite{de1978numerically,GoLo:1996,Gi:TAMS2002}, orthogonal polynomials on the unit circle are a crucial tool for analysis of the continuous one-dimensional Schr\"odinger operator \cite{Gi:SIAMinpress}. The discrete one-dimensional Schr\"odinger operator involves something in between! Indeed, the polynomials (\ref{key-polynomials}) appear not to have been previously studied.  

The hybrid character of equation (\ref{discrete-2}) is manifest in other ways as well.  For instance, its continuous spectrum is a feature shared with its continuous counterpart, while the irrelevance of bound states to inverse scattering, as implied by Theorem~\ref{thm-analytic-inversion}, is a feature common to finite Jacobi operators. Dispersion of travelling waves distinguishes (\ref{discrete-2}) from both its continuous and finite-dimensional counterparts: there is no dispersion in the continuous setting (where $k=\sigma$); and there are no travelling waves in the finite setting.

\subsection{Poles of $\mu_{\oJ}$}
\begin{figure}[htb!]
{
(a)\includegraphics[width=.25\textwidth]{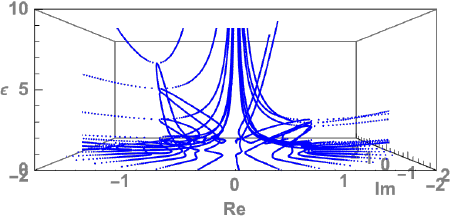}
(b)\includegraphics[width=.25\textwidth]{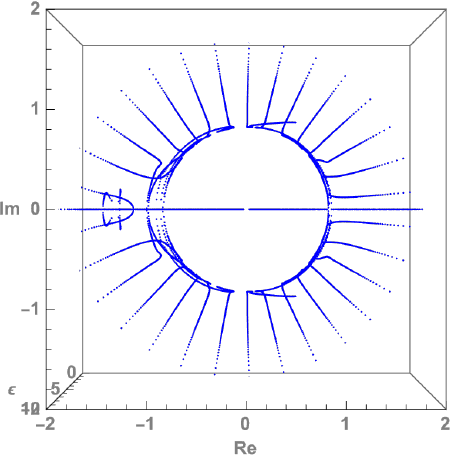}
(c)\includegraphics[width=.25\textwidth]{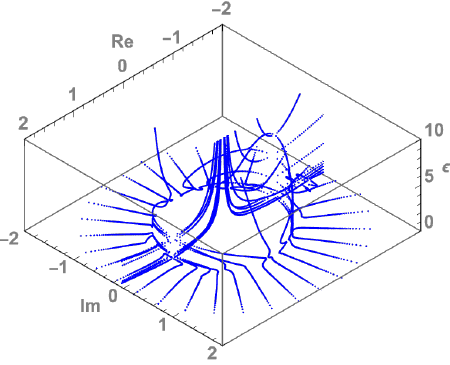}
}
\caption{Various views of the poles of the $\mu$-function for $\oJ_{\varepsilon}$ \eqref{defJeps}, as $\varepsilon$ varies on a grid from $0.01$ to $10$. (a) complex plane orthogonal to the page, (b) looking down on the complex plane, (c) orthographic view. Spacing of the dots corresponds to rate of change of the poles with respect to $\varepsilon$.
}\label{fig-pole-tracking}
\end{figure}

Here we consider, for a fixed perturbed Schr\"odinger operator $\oJ$ as a random perturbation of $\defJz$, 
\begin{equation}
\oJ_\varepsilon\defeq \defJz+\varepsilon\left(\oJ-\defJz\right). 
\label{defJeps}
\end{equation}
The flow of the poles of $\mu_{\oJ_\varepsilon}$ (which are the zeros of $z^LB(z)$ in (\ref{B-1})) as $\varepsilon$ increases exhibits some very interesting behaviour. 
As an illustration, we consider the perturbation length $\nuL=16$. 
The flow of the poles is depicted in Figure~\ref{fig-pole-tracking}. Observe that the poles come from infinity toward the unit circle along approximately equally-spaced radial paths. Behaviour near the unit circle becomes complicated, but they appear to only enter the unit circle along the real axis.  Significant parts of their trajectories are on or very near the unit circle. This intriguing phenomenon calls for further investigation. 

\subsection{Open questions and further research}
A comparison with the continuous and finite settings suggests a number of open questions. 

For example, it is known that in the case of continuous one-dimensional Schr\"odinger operators, even if the bound states and norming constants are given, the transmission coefficient does not determine the potential. And in the finite setting, the analogue of the transmission coefficient gives a lower bound on the dimension of solutions, but no upper bound \cite[Thm.~6.2]{Gi:TAMS2002}.  An open question for (\ref{discrete-2}) is to what extent the transmission coefficient alone determines the perturbation. The formula for $T_1$ and $T_2$ given in Theorem~\ref{thm-S-formula} lays the groundwork for this question, but we have not pursued it further in this article. 

Another open question is to characterize the class of all possible $S$-matrices for equations of the form (\ref{discrete-2}), as has been done in the continuous setting \cite{Fa:1964,DeTr:1979}. 

A classical fact about orthogonal polynomials on the line, such as those associated to finite and semi-infinite Jacobi operators, is that they have simple roots. Our algorithms are predicated on the generic assumption that the polynomials (\ref{key-polynomials}) have simple roots, which is consistently the case in numerical experiments. But it remains an open question whether multiple roots of polynomials of the form (\ref{key-polynomials}) are theoretically possible. 

Other open questions relate to computational aspects of the inverse problem. 
One such question concerns a curious feature of the right-most columns of Figures~\ref{fig:timing_error_comparison} that the largest error does not occur exactly where one would expect. We suspect this is related to complicated behaviour of the poles of the function $\mu_{\oJ}$ at large perturbation scales, whereby they come close to the unit circle as depicted in Figure~\ref{fig-pole-tracking}. 

Another question related to computational aspects concerns the minimal number of values $R_2(k_j)$ needed to recover $\oJ$. This did not arise in present article since we were concerned first to determine $\nuL$ under the assumption that the number of given data points is of a larger order, but the question is nonetheless of interest.  And supposing one is given the minimal number of data points, how does one compute $\oJ$ if the data points are not equally spaced on the circle, say? The multi-point Pad\'e approximation \cite[\S5.5.2]{JoTh:1980} is a natural candidate, but the details have yet to be worked out. 

\section*{Acknowledgements}
This work was initiated during the stay in 2023 of both authors at the Isaac
Newton Institute (INI) for Mathematical Sciences, Cambridge. The authors would
like to thank INI for support and hospitality during the programme ‘Mathematical
theory and applications of multiple wave scattering’ (MWS). This work was supported by EPSRC grant no EP/R014604/1.
Work of both authors, for the same visit to INI, was partially supported by a grant from the Simons Foundation.
Research by PCG was partially supported by NSERC grant DG RGPIN-2022-04547.


\appendix

\section{Analysis of poles\label{sec-poles}}

Let $\ddd\subset\complex$ denote the open unit disk in the complex plane, and let $S^1$ denote its boundary, the unit circle.  Let 
\begin{equation}\label{phifunc}
\varphi:\complex\rightarrow\complex
\end{equation}
be a rational function having real coefficients, i.e., such that $\overline{\varphi(z)}=\varphi\bigl(\bar{z}\bigr)$.
Given $\varphi|_{S^1}$, we seek to determine: (i) the number of poles of $\varphi$; and (ii) the value $\varphi(0)$.  For ease of presentation we consider just the generic case where the poles of $\varphi$ are simple.  In short, the requisite information is encoded in the Fourier coefficients of $\varphi$, which can be computed from $\varphi|_{S^1}$. 

Fix notation as follows.  
Set 
\begin{equation}\label{fourier-coefficients}
c_j=\frac{1}{2\pi}\int_{-\pi}^\pi \varphi(e^{ik})e^{-ijk}\,dk\qquad(j\in\integer).
\end{equation}
If $\varphi$ has no poles in the closed unit disk $\overline{\ddd}$, then
$\varphi$ is harmonic on $\overline{\ddd}$, and so $\varphi(0)=c_0$.  Suppose henceforth that $\varphi$ has $m\geq1$ poles, denoted
\begin{equation}\label{poles}
\lambda_1,\ldots,\lambda_m\quad\mbox{ where }\quad |\lambda_i|\leq |\lambda_j|\quad\mbox{ for every }\quad i<j,
\end{equation}
with precisely $\kappa\geq1$ poles in $\overline{\ddd}$. 
Thus the poles are indexed such that
\[
|\lambda_j|\leq 1\quad\mbox{ for }\quad 1\leq j\leq \kappa,\quad\mbox{ and }\quad |\lambda_j|>1\quad\mbox{ for }\quad \kappa<j\leq m. 
\] 
Note that if at least one pole lies on the unit circle, then the integral (\ref{fourier-coefficients}) is singular, but it is nevertheless well defined in a principal value sense.  The principal value interpretation is to be understood whenever poles of $\varphi$ lie on $S^1$. 
For $i,j\geq 1$, set 
\begin{equation}\label{hankel}
H_{i,j}=\begin{pmatrix}
c_{-1}&c_{-2}&\cdots&c_{-j}\\
c_{-2}&c_{-3}&\cdots&c_{-j-1}\\
&&\vdots&\\
c_{-i}&c_{-i-1}&\cdots&c_{-i-j+1}
\end{pmatrix},
\widetilde{H}_{i,j}=\begin{pmatrix}
c_{-2}&c_{-3}&\cdots&c_{-j-1}\\
c_{-3}&c_{-4}&\cdots&c_{-j-2}\\
&&\vdots&\\
c_{-i-1}&c_{-i-2}&\cdots&c_{-i-j}
\end{pmatrix},
V_j=\begin{pmatrix}1&\lambda_1&\cdots&\lambda_1^{j-1}\\
1&\lambda_2&\cdots&\lambda_2^{j-1}\\
&&\vdots&\\
1&\lambda_\kappa&\cdots&\lambda_\kappa^{j-1}\end{pmatrix}.
\end{equation}
Thus $H_{i,j}$ and $\widetilde{H}_{i,j}$ are $i\times j$ Hankel matrices constructed from the Fourier coefficients of $\varphi$ having negative index, and $V_j$ is a $\kappa\times j$ Vandermonde matrix constructed from the poles inside the closed unit disk.  Observe that $\varphi$ has a unique representation of the form
\begin{equation}\label{partial-fraction}
\varphi=\varphi_a+\varphi_s\quad\mbox{ where }\quad \varphi_s(z)=\sum_{j=1}^\kappa\frac{w_j}{z-\lambda_j}, \quad 0\neq w_j\in\real\quad (1\leq j\leq \kappa),
\end{equation}
and where $\varphi_a$ is analytic on the closed unit disk.  It follows that for negative indices $j$,
\begin{equation}\label{negative-indices}
c_j=\frac{1}{2\pi}\int_{-\pi}^\pi \bigl(\varphi_a(e^{ik})+\varphi_s(e^{ik})\bigr)e^{-ijk}\,dk
=\frac{1}{2\pi}\int_{-\pi}^\pi \varphi_s(e^{ik})e^{-ijk}\,dk\qquad(j<0).
\end{equation}
Note
\begin{equation}\label{constant-term}
\frac{1}{2\pi}\int\limits_{-\pi}^\pi\frac{1}{e^{ik}-\lambda_m}\,dk=\left\{
\begin{array}{cc}
0&\mbox{ if }|\lambda_m|<1\\[5pt]
-\displaystyle\frac{1}{2\lambda_m}&\mbox{ if }|\lambda_m|=1\\[10pt]
-\displaystyle\frac{1}{\lambda_m}&\mbox{ if }|\lambda_m|>1
\end{array}\right.
\end{equation}
\begin{equation}\label{negative-index}
\mbox{ and for $j\geq 1$,}\qquad
\frac{1}{2\pi}\int\limits_{-\pi}^\pi\frac{e^{ijk}}{e^{ik}-\lambda_m}\,dk=\left\{
\begin{array}{cc}
\lambda_m^{j-1}&\mbox{ if }|\lambda_m|<1\\[5pt]
\frac{1}{2}\lambda_m^{j-1}&\mbox{ if }|\lambda_m|=1\\[5pt]
0&\mbox{ if }|\lambda_m|>1
\end{array}\right..\rule{54pt}{0pt}
\end{equation}
Referring to (\ref{partial-fraction}), set 
\begin{equation}\label{diagonal-matrices}
W=\diag(\varepsilon_1w_1,\ldots,\varepsilon_\kappa w_\kappa)\quad\mbox{ and }\quad \Lambda=\diag(\lambda_1,\ldots,\lambda_\kappa),
\end{equation}
where 
\begin{equation}\label{defeps}
\varepsilon_i=1\text{ if }|\lambda_i|<1\text{ and }\varepsilon_i=1/2\text{ if }|\lambda_i|=1\quad(1\leq i\leq \kappa).
\end{equation}
It follows from (\ref{partial-fraction}), (\ref{negative-indices}) and (\ref{negative-index}) that for every $m,j\geq 1$, 
\begin{equation}\label{hankel-decomposition}
H_{m,j}=V_m^\top WV_j\quad\mbox{ and }\quad \widetilde{H}_{m,j}=V_m^\top W\Lambda V_j. 
\end{equation}
Observe by (\ref{hankel-decomposition}) that 
\begin{equation}\label{rank}
\rank H_{m,j}=\min\{m,j,\kappa\}\qquad(m,j\geq 1). 
\end{equation}
The matrices $H_{m,j}$ and $\widetilde{H}_{m,j}$ can be computed from $\varphi_{S^1}$, since they are built from the Fourier coefficients of $\varphi$.  Moreover, they determine the number of poles of $\varphi$ inside $\overline{\ddd}$ and the value $\varphi(0)$, as follows. 
\begin{prop}\label{prop-poles}
Let $\varphi$ be a rational function having simple poles, $\kappa\geq1$ of which are in $\overline{\ddd}$.  Define $H_{m,j}$, $\widetilde{H}_{m,j}$, $V_j$ and $W$ as in (\ref{hankel})
and (\ref{diagonal-matrices}).  If $m,j\geq \kappa$, then 
\begin{equation}\label{kappa-determination}
\kappa=\rank H_{m,j}.
\end{equation}
Set 
\begin{equation}\label{spec-K}
{K}\defeq\widetilde{H}_{\kappa,\kappa}\bigl(H_{\kappa,\kappa}\bigr)^{-1}. 
\end{equation}
Then
\begin{equation}\label{spec-K2}
\{\lambda_1,\ldots,\lambda_\kappa\}=\spec {K}\thickspace\mbox{ and }\thickspace W=\diag(\varepsilon_1w_1,\ldots,\varepsilon_\kappa w_\kappa)=(V_\kappa^\top)^{-1}H_{\kappa,\kappa}V_\kappa^{-1}.
\end{equation}
The value $\varphi(0)$ is given by the formula
\begin{equation}\label{varphi-0}
\varphi(0)=\frac{1}{2\pi}\int\limits_{-\pi}^\pi\varphi(e^{ik})\,dk-\sum_{j=1}^\kappa\frac{\varepsilon_jw_j}{\lambda_j}.
\end{equation}
\end{prop}
For a rational function $\varphi$ having real coefficients, the number of its poles outside $\overline{\ddd}$ can be computed as follows. The function $\bar\varphi|_{S^1}$ has a unique rational extension to $\complex$; denote this extension $\varphi^\ast$.  Then poles of $\varphi^\ast$ in $\ddd$ correspond to poles of $\varphi$ outside $\overline{\ddd}$.  Thus the total number of poles of $\varphi$ is the number of poles of $\varphi$ in $\overline{\ddd}$ plus the number of poles of $\varphi^\ast$ in $\ddd$, both of which can be computed by dint of Proposition~\ref{prop-poles}. 

\begin{remark}\label{remark-poles}
It follows from the definition of ${K}$ in Proposition~\ref{prop-poles} that ${K}H_{\kappa,j}=\widetilde{H}_{\kappa,j}$ for any $j\geq1$.  Thus the Fourier coefficients of $\varphi$ corresponding to negative indices satisfy the linear recurrence
\begin{equation}\label{linear-recurrence}
\begin{pmatrix}c_{-j-1}\\ \vdots\\ c_{-j-\kappa}\end{pmatrix}={K}\begin{pmatrix}c_{-j}\\ \vdots\\ c_{-j-\kappa+1}\end{pmatrix}\qquad(j\geq 1)
\end{equation}
whereby the matrix ${K}$ acts as a right shift operator on blocks of length $\kappa$. This fact is often attributed to Prony  \cite{Pr:1795}. 
\end{remark}

\section{Algorithms\label{sec-algorithms}}

Proposition~\ref{prop-poles} provides the theoretical basis for evaluation of $\varphi(0)$ given $\varphi|_{S^1}$. Practical implementation is facilitated by taking into account 
numerical considerations, as follows. 

In practice, the rank $\kappa$ as specified in (\ref{kappa-determination}) is obtained by counting the singular values of $H_{m,j}$ larger than a specified tolerance. 

Consider the matrix $K$ defined In Proposition~\ref{prop-poles} by the formula \eqref{spec-K} $K=\widetilde{H}_{\kappa,\kappa}\bigl(H_{\kappa,\kappa}\bigr)^{-1}$. For a non-zero vector $y=H_{\kappa,\kappa}x\in\complex^\kappa$ and a scalar $\lambda\in\complex$,
\begin{equation}\label{generalized-eigenvalue}
Ky=\lambda y\quad\Leftrightarrow\quad \widetilde{H}_{\kappa,\kappa}x=\lambda H_{\kappa,\kappa}x.
\end{equation}
Thus eigenvalues of $K$ are generalized eigenvalues of the matrix pair $\bigl(\widetilde{H}_{\kappa,\kappa}, H_{\kappa,\kappa}\bigr)$.  The latter formulation is computationally advantageous in determining $\spec K$ since, for instance, it does not require computation of the inverse matrix $H_{\kappa,\kappa}^{-1}$; see \cite[\S7.7]{GoLo:1996}. We use the notation $\spec(A,B)$ for the vector of generalized eigenvalues of the matrix pair $(A,B)$; see line~\ref{spec-line} in Algorithm~\ref{alg-pole-correction}.

Regarding the right-hand equation in (\ref{spec-K}), observe that the left-most column of the equivalent equation
\[
H_{\kappa,\kappa}=V_\kappa^\top WV_\kappa
\]
yields
\begin{equation}\label{W-alt}
\begin{pmatrix}c_{-1}\\ c_{-2}\\ \vdots\\ c_{-\kappa}\end{pmatrix}
=V_\kappa^\top\begin{pmatrix}\varepsilon_1w_1\\ \varepsilon_2w_2\\ \vdots\\ \varepsilon_\kappa w_\kappa\end{pmatrix}.
\end{equation}
Equation (\ref{W-alt}) enables efficient computation of $W$, given $\spec K$ and the Fourier coefficients of $\varphi$, $c_{-j}$ $(1\leq j\leq \kappa)$. 

Lastly, note that equation (\ref{kappa-determination}) in Proposition~\ref{prop-poles} counts poles inside the closed unit disk.  The number of poles lying on the unit circle $S^1$ (within a given tolerance) can be estimated by computation of $\spec K$.  In order to count the total number of poles of $\varphi$, then, one needs to determine the number of poles of $\varphi$ in $\overline{\ddd}$ and the number of poles of $\varphi^\ast$ in $\overline{\ddd}$, and then subtract the number of poles of $\varphi$ on $S^1$, which are double counted, since they are poles of both $\varphi$ and $\varphi^\ast$.   Computation of $\spec K$ can be avoided if one is content to obtain an upper estimate on the number of poles in $\varphi$ by adding together the numbers of poles of $\varphi$ and $\varphi^\ast$ in $\overline{\ddd}$.  This is in fact preferable because it avoids the issue of deciding when a number $\lambda\in\spec K$ should be considered to have modulus 1 up to numerical error. 

\floatname{algorithm}{Algorithm}
\renewcommand{\thealgorithm}{1} 
\begin{algorithm}[htb!]
\caption{Pole count inside the closed unit disk}\label{alg-pole-count}
\begin{algorithmic}[1]
\Procedure{PoleCount}{$v,N,\epsilon$}
\State {\bf Input} $(v,N,\epsilon)\in\complex^d\times\{1,\ldots,\lfloor (d+1)/2\rfloor\}\times\real_+$
\State $\omega\gets e^{2\pi i/d}$
\State $\Omega(m,j)\gets\omega^{mj}\quad(1\leq m\leq 2N-1;\thickspace 1\leq j\leq d)$
\State $\gamma\gets -\frac{e^{-i\pi/d}}{d}\Omega v$
\State $H(m,j)\gets\gamma_{m+j-1}\quad(1\leq m,j\leq N)$  
\State $\sigma\gets$ vector of singular values of $H$ such that $\sigma_1\geq\cdots\geq\sigma_N\geq 0$
\If{$\sigma_1\leq\epsilon$}
\State $\bigl(\kappa,\widetilde{H},H,\gamma\bigr)\gets(0,\emptyset,\emptyset,\emptyset)$
\Else
\State $\kappa\gets\max\{j\,|\,\sigma_j>\epsilon\}$
\State $\widetilde{H}(m,j)\gets\gamma_{m+j}\quad(1\leq m,j\leq \kappa)$
\State $H\gets$ $\kappa\times\kappa$ leading principal submatrix of $H$
\State $\gamma\gets(\gamma_1,\ldots,\gamma_\kappa)^\top$
\EndIf
\State {\bf Output} $\bigl(\kappa,\widetilde{H},H,\gamma\bigr)\in\integer_+\times\complex^{\kappa\times\kappa}\times\complex^{\kappa\times\kappa}\times\complex^\kappa$
\EndProcedure
\end{algorithmic}
\end{algorithm}
Note that the output of Algorithm~\ref{alg-pole-count} consists not just of the pole count $\kappa$, but also matrices $\widetilde{H}, H$ and the vector $\gamma$; the latter items are retained to reduce computation when Algorithm~\ref{alg-pole-count} is called by Algorithm~\ref{alg-pole-correction}, below.

\floatname{algorithm}{Algorithm}
\renewcommand{\thealgorithm}{2} 
\begin{algorithm}[htb!]
\caption{Pole correction for a rational function sampled on $S^1$}\label{alg-pole-correction}
\begin{algorithmic}[1]
\Procedure{PoleCorrection}{$v,N,\epsilon$}
\State {\bf Input} $(v,N,\epsilon)\in\complex^d\times\{1,\ldots,\lfloor (d+1)/2\rfloor\}\times\real_+$
\State $\bigl(\kappa,\widetilde{H},H,\gamma\bigr)\gets\mbox{\sc PoleCount}(v,N\epsilon)$
\If{$\kappa=0$}
\State $C\gets 0$ 
\Else
\State $\lambda\gets \spec\bigl(\widetilde{H},H\bigr)$ \label{spec-line}
\State $V(m,j)\gets\lambda_m^{j-1}\quad(1\leq m,j\leq \kappa)$
\State $y\gets \solve V^\top y=\gamma$
\State $C\gets-\sum_{j=1}^\kappa\frac{y_j}{\lambda_j}$
\EndIf
\State {\bf Output} $C\in\complex$
\EndProcedure
\end{algorithmic}
\end{algorithm}

\end{document}